\documentclass[11pt]{amsart}

\usepackage{amssymb}
\usepackage{amscd}

\usepackage[all]{xy}
\begin{document}
\title{Quasi-lines and their degenerations}
\author{Laurent BONAVERO, Andreas H\"ORING}
\date{February 16th, 2007}
\maketitle
\noindent
\def\restriction{\string |}
\newcommand{\pp}{\rm ppcm}
\newcommand{\pg}{\rm pgcd}
\newcommand{\Ker}{\rm Ker}
\newcommand{\C}{{\mathbb C}}
\newcommand{\Q}{{\mathbb Q}}
\newcommand{\GL}{\rm GL}
\newcommand{\SL}{\rm SL}
\newcommand{\diag}{\rm diag}
\newcommand{\N}{{\mathbb N}}
\def\finpreuve
{\hskip 3pt \vrule height6pt width6pt depth 0pt}

\newtheorem{lemma1}{}[section]
\newtheorem*{theonn}{Theorem}

\newenvironment{lemm}{\begin{lemma1}{\bf Lemma.}}{\end{lemma1}}
\newenvironment{example}{\begin{lemma1}{\bf Example.}\rm}{\end{lemma1}}
\newenvironment{examples}{\begin{lemma1}{\bf Examples.}\rm}{\end{lemma1}}
\newenvironment{question}{\begin{lemma1}{\bf Question.}\rm}{\end{lemma1}}
\newenvironment{abs}{\begin{lemma1}\rm}{\end{lemma1}}
\newenvironment{theo}{\begin{lemma1}{\bf Theorem.}}{\end{lemma1}}
\newenvironment{prop}{\begin{lemma1}{\bf Proposition.}}{\end{lemma1}}
\newenvironment{cor}{\begin{lemma1}{\bf Corollary.}}{\end{lemma1}}
\newenvironment{rem}{\begin{lemma1}{\bf Remark.}\rm}{\end{lemma1}}
\newenvironment{defi}{\begin{lemma1}{\bf Definition.}}{\end{lemma1}}
\newenvironment{construction}{\begin{lemma1}{\bf Construction.}}{\end{lemma1}}
\newenvironment{conjecture}{\begin {lemma1}{\bf Conjecture.}}{\end{lemma1}}
\newenvironment{problem}{\begin{lemma1}{\bf Problem.}}{\end{lemma1}}

\newcommand{\CC}{{\mathbb C}}
\newcommand{\ZZ}{{\mathbb Z}}
\newcommand{\RR}{{\mathbb R}}
\newcommand{\QQ}{{\mathbb Q}}
\newcommand{\FF}{{\mathbb F}}
\newcommand{\PP}{{\mathbb P}}
\newcommand{\NN}{{\mathbb N}}
\newcommand{\codim}{\operatorname{codim}}
\newcommand{\Ho}{\operatorname{Hom}}
\newcommand{\Pic}{\operatorname{Pic}}
\newcommand{\NE}{\operatorname{NE}}
\newcommand{\Nun}{\operatorname{N}}
\newcommand{\card}{\operatorname{card}}
\newcommand{\Hilb}{\operatorname{Hilb}}
\newcommand{\mult}{\operatorname{mult}}
\newcommand{\vol}{\operatorname{vol}}
\newcommand\sO{{\mathcal O}}
\newcommand{\divi}{\operatorname{div}}
\newcommand{\pr}{\operatorname{pr}}
\newcommand{\con}{\operatorname{cont}}
\newcommand{\ima}{\operatorname{Im}}
\newcommand{\rk}{\operatorname{rk}}
\newcommand{\Exc}{\operatorname{Exc}}
\newcounter{subsub}[subsection]
\def\thesubsub{\thesubsection .\arabic{subsub}}
\def\subsub#1{\addtocounter{subsub}{1}\par\vspace{3mm}
\noindent{\bf \thesubsub ~ #1 }\par\vspace{2mm}}
\def\coker{\mathop{\rm coker}\nolimits}
\def\pr{\mathop{\rm pr}\nolimits}
\def\im{\mathop{\rm Im}\nolimits}
\def\hfl#1#2{\smash{\mathop{\hbox to 12mm{\rightarrowfill}}
\limits^{\scriptstyle#1}_{\scriptstyle#2}}}
\def\vfl#1#2{\llap{$\scriptstyle #1$}\big\downarrow
\big\uparrow
\rlap{$\scriptstyle #2$}}
\def\diagram#1{\def\normalbaselines{\baselineskip=0pt
\lineskip=10pt\lineskiplimit=1pt}   \matrix{#1}}
\def\limind{\mathop{\oalign{lim\cr
\hidewidth$\longrightarrow$\hidewidth\cr}}}

\long\def\InsertFig#1 #2 #3 #4\EndFig{
\hbox{\hskip #1 mm$\vbox to #2 mm{\vfil\includegraphics{#3}}#4$}}
\long\def\LabelTeX#1 #2 #3\ELTX{\rlap{\kern#1mm\raise#2mm\hbox{#3}}}


\newcommand{\merom}[3]{\ensuremath{#1:#2 \dashrightarrow #3}}
\newcommand{\holom}[3]{\ensuremath{#1:#2  \rightarrow #3}}
\newcommand{\fibre}[2]{\ensuremath{#1^{-1} (#2)}}
\newcommand{\Z}{{\mathbb Z}}
\newcommand{\chow}[1]{\ensuremath{\mathcal{C}(#1)}}
\newcommand{\ns}{\operatorname{ns}}
\newcommand{\red}{\operatorname{red}}

\newcommand\sI{{\mathcal I}}

\setcounter{tocdepth}{1}
\renewcommand\thesubsection{\thesection.\Alph{subsection}}


{\let\thefootnote\relax
\footnote{
\textbf{A.M.S.~classification:} 14E30, 14J10, 14J30, 14J40, 14J45. 
}}

\vspace{-1cm}


\begin{center}
\begin{minipage}{130mm}
\scriptsize

{\bf Abstract.} In this paper we study the structure
of manifolds that contain a quasi-line and give some evidence 
towards the fact that the irreducible components of degenerations 
of the quasi-line should determine the Mori cone. We show that the minimality
with respect to a quasi-line yields strong restrictions on fibre space structures
of the manifold.
\end{minipage}
\end{center}

\tableofcontents

\section{Introduction} 

Let $X$ be a complex quasiprojective manifold of dimension $n$.
A quasi-line $l$ in $X$ is a {\it smooth} rational curve $f~: \PP^1 
\hookrightarrow X$
such that $f^*T_X$ is the same as for a line in $\PP^n$, {\em i.e.} 
is isomorphic to 
$${\mathcal O}_{\PP^1}(2) \oplus {\mathcal O}_{\PP^1}(1)^{\oplus n-1}.$$
Although the terminology suggests that quasi-lines are very special objects,
we will see that they appear in a lot of situations.

\medskip

\begin{examples} 
\label{examplesstart}
\begin{enumerate}
\item[(1)] If $X$ is a smooth Fano threefold of index $2$
with ${\rm Pic}(X) = \ZZ H$, where $H$ is very ample,
then a general conic $C$ ({\em i.e.} a curve satisfying $H \cdot C =2$)
is a quasi-line (Oxbury \cite{Ox94}, see also B\u adescu, 
Beltrametti and Ionescu \cite{BBI00}). 
\item[(2)] If $X$ is rationally connected, then there exists
a sequence $X' \to X$ of blow-ups along smooth
codimension $2$ centres such that $X'$ contains a quasi-line (Ionescu and Naie \cite{IN03}).
\item[(3)] If $X$ contains a quasi-line $l$,
let $\pi~: X' \to X$ be a blow-up of $X$ along a smooth subvariety $Z$. 
A general deformation of $l$ does not meet $Z$, 
so it identifies to a quasi-line
in $X'$. 
\end{enumerate}
\end{examples}

Deformation theory shows that rational curves 
whose deformations passing through a fixed point
dominate the manifold have $-K_X$-degree at least $\dim X+1$. Quasi-lines 
can be seen as the rational curves realising the boundary case $-K_X \cdot l=\dim X+1$,
so it is reasonable to ask if the existence of a quasi-line has any
implications on the global structure of the manifold. Yet as the example \ref{examplesstart}(2) shows
this implication can't be much stronger than the rational connectedness 
of $X$, so we will have to make extra restrictions. 
The theory of K\"ahlerian twistor spaces which provided the first motivation for the study of quasi-lines
suggests that the most important class to study are Fano manifolds \cite{Hi81}.
The cone theorem then shows that the Mori cone ${\NE}(X)$
of a Fano manifold is closed and polyhedral, the extremal rays being generated by classes of rational curves.
If the Picard number is at least 2, the class of a quasi-line does not generate an extremal ray,
but we have the natural following question.

\medskip

\begin{question}
\label{questionextremalrays} For a Fano manifold $X$ containing a quasi-line $l$,
do the numerical classes of irreducible components of degenerations 
of $l$ generate the Mori cone ${\NE}(X)$~? 
\end{question}

\medskip

In this question, it is hopeless to expect that the numerical classes 
of irreducible components of {\em a single} degeneration 
of $l$ generate the Mori cone ${\NE}(X)$.  
Before asking such a question, one should better verify that degenerations exists.
This is guaranteed by the characterisation of the projective space
by Cho, Miyaoka and Shepherd-Barron, which we restate here in the language
of quasi-lines (see section \ref{sectionbasics} for the terminology).
 
\begin{theo}  \cite{CMS02,Ke01}
\label{theoprojectivespace}
Let $X$ be a projective manifold that contains a quasi-line $l$.
Suppose that for a general point $x$ of $X$, the deformations of $l$ passing through $x$
form an unsplit family of rational curves.
Then $X \simeq \PP^n$ and $l \subset \PP^n$ is a line. 
\end{theo}

First evidence for an affirmative answer to the question comes from the following situation:
suppose that there exists an effective divisor $D \subset X$ such that $D \cdot l = 0$.
Then there exists a degeneration that has an irreducible component in $D$, moreover there exists a birational
Mori contraction whose locus is contained in $D$ 
(see lemma \ref{lemmabirationalcontractions}).
This observation leads us to recall the related notion 
of minimality with respect to a quasi-line, introduced by 
Ionescu {\em et al.} in \cite{BBI00, IV03}.

\begin{defi}
\label{definitionminimal}
Let $X$ be a projective normal $\Q$-factorial (Fano) variety $X$ 
that contains a quasi-line $l \subset X_{\ns}$ (we say that the couple $(X,l)$ is a (Fano) model). 
The variety $X$ is minimal with respect to $l$
if for every effective Cartier divisor $\ D \subset X$, we have $D \cdot l>0$.
This also means that the numerical class of $l$ belongs to the interior
of the cone generated by the classes of moving curves.

Two models $(X,l)$ and $(X',l')$ are equivalent if there are 
Zariski open subsets $U \subset X$ and $U' \subset X'$
containing respectively $l$ and $l'$ and an isomorphism
$\mu~: U \to U'$ such that $\mu(l)=l'$.
\end{defi}

The notions of minimality and equivalent models were introduced 
in order to avoid situations like in example \ref{examplesstart}(3), where $X'$ is clearly not minimal. 
Our first main result is an inverse statement for smooth Fano models of dimension three.

\medskip

\noindent {\bf Theorem~A.}
{\em Let $(X,l)$ be a smooth Fano model of dimension three. 
Then there exists a birational morphism 
\holom{\mu}{X}{X'} onto a $\Q$-factorial projective threefold 
with at most terminal singularities 
such that $\mu$ is an equivalence of models $(X,l) \simeq (X',\mu(l))$, 
and $X'$ is minimal with respect to $\mu(l)$.}

\medskip

The core of this paper considers smooth models $(X,l)$ such that the manifold $X$ 
admits a fibration \holom{\varphi}{X}{Y}. The first task is to characterise
situations where this induces a morphism of models $(X,l) \rightarrow (Y, \varphi(l))$, {\em i.e.} the image 
$\varphi(l) \subset Y_{\ns}$ is a quasi-line.
Lemma \ref{lemmafibrations} shows 
that the study of morphisms of models should be started 
by considering fibrations 
\holom{\varphi}{X}{\PP^1} such that $l$ is a section. 
We will investigate this in subsection \ref{subsectionfibrationcurve} and discover important
differences between models 
with $e(X,l)=1$ (that is there exists exactly one deformation of $l$
through two general points of $X$, see definition \ref{definitioneinvariant}) and those with $e(X,l)>1$.

\medskip

\noindent {\bf Theorem~B.}
{\em
Let $(X,l)$ be a smooth Fano model 
of dimension at most four, and let $\varphi: X \to \PP^1$ be a fibration 
such that $l$ is a section of $\varphi$.
Then the general fibre $F$ contains a quasi-line.
If furthermore $e(X,l)=1$, the variety $X$ is not minimal with respect to $l$.
}

\medskip

If we accept the general idea that the extremal contractions of a Fano model $(X,l)$ are
related to the degenerations of $l$, the minimality with respect to $l$ should
also yield some restriction on the structure of the morphisms of $X$.
An example for a birational morphism 
is the  following characterisation of the projective space
due to B\u adescu, Beltrametti and Ionescu.

\begin{theo} \cite[Thm.4.4, Cor. 4.6]{BBI00}
\label{theoremalmostline}
Let $(X,l)$ be a smooth model such that $X$ is minimal with respect to $l$.
Suppose that there exists a big and nef divisor 
such that $H \cdot l = 1$. Then $X \simeq \PP^n$ and $l \subset \PP^n$ is a line.
\end{theo}

We will show a similar statement in a relative setting. 

\medskip

\noindent {\bf Theorem~C.}
{\em
Let \holom{\varphi}{(X,l)}{(Y,\varphi(l))} be a morphism of smooth models,
and suppose that $\varphi$ is flat of relative dimension $1$.
Assume furthermore that:
\begin{enumerate}
\item[(1)] $X$ is minimal with respect to $l$ and $e(X,l)=1$,
\item[(2)] there exists a big and nef line bundle $H$ on $Y$ 
such that $\varphi^* H \cdot l=1$. 
\end{enumerate}
Then $Y\simeq \PP^{n-1}$ and $X \simeq \PP(E)$ where $E$ is a stable rank two vector bundle over $Y$.
If furthermore $X$ is Fano, we have $\dim X=3$ and $X \simeq \PP(T_{\PP^2})$.
}

\medskip

This result follows from two intermediate results
(theorem \ref{theoremnotconicbundle} and proposition \ref{propositionstability}) which are 
interesting in their own right and classification results of Fano bundles (corollary \ref{corollaryfanobundlem3}
and proposition \ref{propositionp2fanobundle}).

\medskip 

{\bf Acknowledgements.} 
We heartly thank Cinzia Casagrande and Rita Pardini for valuable discussions
and explaining to us the very interesting examples \ref{examplespardini}.
We also thank St\'ephane Druel for his 
proofreading and helping us to repair an error
in the proof of proposition \ref{propositionsmoothing}.
{\em This work has been partially supported by the 3AGC project
of the A.N.R.}

\section{Notation and basic results}
\label{sectionbasics}

We work over the complex field $\C$, topological notions always refer to the Zariski topology.
A variety is an integral scheme of finite type over $\C$, a manifold is a smooth variety.
A fibration is a surjective 
morphism \holom{\varphi}{X}{Y} between normal varieties such that
$\dim X>\dim Y$ and $\varphi_* \sO_X \simeq \sO_Y$,
that is all the fibres are connected. Fibres are always scheme-theoretic fibres.
For general definitions we refer to Hartshorne's book \cite{Ha77}, we will also use the standard terminology
of Mori theory and deformation theory as explained in \cite{Deb01, Ko96}.

\subsection{Deformations and degenerations}
Let $X$ be a projective variety, and let 
$C \subset X$ be an integral projective curve.
Identify 
$C$ to its fundamental 
cycle\footnote{Throughout the whole paper, We will not distinguish between an effective cycle and its support.}, and 
suppose that the Chow scheme $\chow{X}$ is irreducible at 
the point $[C]$. Denote by ${\mathcal H}$ the normalisation of the 
unique irreducible 
component of $\chow{X}$ containing $[C]$. 
Furthermore we have the incidence variety $\mathcal U \subset \mathcal H \times X$ endowed with 
two natural morphisms 
\holom{q}{\mathcal U}{\mathcal H} and \holom{p}{\mathcal U}{X}, in particular  
$q$ is equidimensional. 

\begin{defi}
A deformation $C'$ of $C$ is an integral curve such that $[C'] \in \mathcal H$. 
Let $h \in \mathcal H$ be a point 
that is not a deformation of $C$, {\em i.e.}
the corresponding cycle $C_h$ is reducible or non-reduced, then $C_h$ will be called 
a degeneration of $C$. Denote by $\mathcal H^* \subset \mathcal H$ 
the open subset parametrizing deformations of $C$.  
The family $\mathcal H$ is said to be unsplit if $\mathcal H^* = \mathcal H$.

Let $Z \subset X$ be a subvariety. We say that the deformations of $C$ dominate (resp. cover) $Z$ if
$p(\fibre{q}{\mathcal H^*}) \cap Z$ is dense in $Z$ (resp. contains $Z$).
\end{defi}

Let $(X,l)$ be a smooth model, and let $x \in l \subset X$ a point.
Since the normal bundle of $l$ is ample, the 
Chow scheme $\chow{X}$ is smooth of dimension $2 \dim X-2$ at $[l]$. 
Therefore there exists
a unique irreducible component ${\mathcal H}$ of $\chow{X}$ 
containing $[l]$. 
The subscheme $q(\fibre{p}{x}) \subset \mathcal H$ is 
smooth at $[l]$ and of dimension $n-1$. We denote by
${\mathcal H}_x$ the normalisation of the irreducible component of $q(\fibre{p}{x})$ 
that contains $[l]$. We denote by $\mathcal U_x$ the normalisation of the corresponding incidence variety
in ${\mathcal H}_x \times X$.
In order to ease the notation we denote the restriction of $p$ and $q$ to $\mathcal U_x$ by the same letter.

\begin{equation}
\label{basicdiagram}
 \xymatrix{ 
{\mathcal U}_x \ar[d]^{q}\ar[r]^{p} & X \\
{\mathcal H}_x & } 
\end{equation}
 
The general fibre of $q$ is a smooth $\PP^1$
and its image in $X$ is still a quasi-line. 
Let $h \in \mathcal H_x$ be a point that parametrises a degeneration of $l$,
then we'll denote the corresponding $1$-cycle by $\sum_i \alpha_i l_i$
or $l \sim \sum_i \alpha_i l_i$, and 
call it a degeneration of $l$ with fixed point $x$. 

One special feature of quasi-lines is that given a quasi-line and two general points, 
there exists finitely many deformations of the quasi-line passing through these two points.
This fact is formalised in the next notion, introduced by Ionescu 
and Voica.

\begin{defi}
\label{definitioneinvariant}
Let $(X,l)$ be a model, and let $x \in X$ be a general point. Let \holom{p}{{\mathcal U}_x}{X}
be the morphism from diagram \ref{basicdiagram}, then we define
\[
e(X,l):= \deg(p).
\]
\end{defi}

Since the evaluation morphism $p$ is surjective and generically finite, 
the invariant $e$ is well-defined.

\subsection{A lemma}
The following technical lemma shows that if $e(X,l)=1$, the minimality with respect to $l$ yields
a strong restriction on the degenerations of $l$.

\begin{lemm}
\label{lemmabirationalevaluation}
Let $(X,l)$ be a smooth model such that $e(X,l)=1$.
Let $x \in X$ be a general point, and 
let $D \subset X$ be a prime effective divisor  
such that
$x \notin D$. Suppose that through a general point of $D$ passes an irreducible component
$l_i$ of a degeneration $\sum_i \alpha_i l_i$ 
of $l$ with fixed point $x$ such that $l_i \subset D$.
Then $D \cdot l = 0$, in particular $X$ is not minimal with respect to $l$.
\end{lemm}

{\em Proof.} We use the notation of the basic diagram \ref{basicdiagram}.
By definition $e(X,l)=1$ implies that 
the map $p:~ {\mathcal U}_x \to X$
is birational.
We argue by contradiction and suppose that $D \cdot l > 0$. 
This implies that for a general point $y \in D$, there exists a 
quasi-line through $x$ and $y$. 
By hypothesis,
there exists also a degeneration of $l$ through $y$. In particular
\fibre{p}{y} is not a singleton, so $y$ is contained in the image of the
exceptional locus of $p$.
This is impossible since this image has codimension at least 2.
\finpreuve

\section{Minimality and birational contractions}

Question \ref{questionextremalrays} asks if there are is a link between a quasi-line $l \subset X$
and the extremal contractions of $X$. If $X$ is minimal with respect to $l$, 
section \ref{sectionfibrations} will provide some evidence in special cases. 
If $X$ is not minimal, we can make a surprisingly simple observation formulated in the next lemma. Recall first that 
an elementary (Mori extremal) contraction $\varphi~: X \to Y$
from a Fano manifold $X$ to a normal variety $Y$ is a map with connected
fibres which contracts exactly the curves whose numerical class belong
to a given extremal ray of the Mori cone ${\rm NE}(X)$ of $X$.

\begin{lemm}
\label{lemmabirationalcontractions}
Let $(X,l)$ be a Fano model, and let $D$ be an effective prime divisor on $X$ such that
$D \cdot l=0$. Then there exists a Mori contraction of birational type
whose exceptional locus is contained in $D$. 
In particular, if all the elementary contractions of $X$ are of fibre
type, then $X$ is minimal with respect to $l$.
\end{lemm}

{\em Proof.} Fix a general point $x \in X \setminus D$,
then the deformations and 
degenerations of $l$ with fixed point $x$ cover $X$. 
Since $D \cdot l=0$, no deformation of $l$ meets $D$.
For any point $y \in D$, we thus can write
$l \sim \sum_{i=1}^k \alpha_i l_i$ with $\alpha_i \in \N ^*$ 
such that one component passes through $y$. Since the cycle 
also passes through $x$,
there exists a component $l_i$ that 
meets $D$ and is not contained in it, so $D \cdot l_i>0$.
Therefore there exists another component, say $l_1$, that satisfies
$D \cdot l_1 <0$. Since $X$ is Fano we have a decomposition
\[
l_1 = \sum_j \beta_j \Gamma_j,
\]
where the $\Gamma_j$ are generators of the extremal rays of the Mori cone
${\rm NE}(X)$ and $\beta_j > 0$.
Since $D \cdot l_1 <0$  there exists at least one
extremal rational curve $\Gamma_j$ such that $D \cdot \Gamma_j <0$. The 
corresponding extremal contraction has its exceptional locus contained
in $D$, therefore is birational.
\finpreuve 

\medskip

\begin{example}
\label{examplesurfaces}
The only smooth model $(S,l)$ of dimension two
that is minimal with respect to $l$ is
the model $(\PP^2, line)$. Indeed consider the deformations of $l$ passing through
a general point $x$. If $(S,l)$ is not $(\PP^2, line)$, 
there exists 
by (the much easier $2$-dimensional version of) 
theorem \ref{theoprojectivespace} 
at least one
degeneration $\sum_{i=1}^k \alpha_i l_i$. Hence
\[
1 = l \cdot l = l \cdot \sum_{i=1}^k \alpha_i l_i
\]
implies that $l \cdot l_i=0$ for some $i$.
More generally if $(S,l)$ is a smooth model of dimension two,
then $S$ is the blow-up of $\PP^2$ at finitely many
points, $l$ being the
inverse image of a general line of $\PP^2$ \cite[Prop.1.21]{IV03}. 
\end{example}

\medskip

In higher dimension 
the structure of the birational contractions \holom{\mu}{X}{X'} becomes more 
complicated. In particular the contraction might be small (which means that
its exceptional locus in $X$ might be of codimension bigger than two). 
In this case $\mu(D) \subset X'$ would be a Weil divisor such that $\mu(l)$ does not meet $\mu(D)$, thus
contracting would not improve the situation and it is not clear if there exist an equivalent model
which is minimal with respect to the quasi-line. 
The next theorem shows that everything works well in dimension three.
   
\begin{theo}
\label{theoremminimalmodel}
Let $(X,l)$ be a smooth Fano model of dimension three. 
Then there exists a birational morphism 
\holom{\mu}{X}{X'} on a $\Q$-factorial projective threefold 
with at most terminal singularities 
such that $\mu$ is an equivalence of models $(X,l) \simeq (X',\mu(l))$, 
and $X'$ is minimal with respect to $\mu(l)$.
\end{theo}  

{\em Proof.} 
Let $D_1, \ldots, D_k$ be the pairwise different effective 
prime divisors such that $D_i \cdot l = 0$. 
Then by lemma \ref{lemmabirationalcontractions}, 
there exists for every $i \in \{ 1, \ldots, k \}$ an elementary contraction 
of birational type
\holom{\mu_i}{X}{X_i} such that the exceptional 
locus is contained in $D_i$. 
Since birational contractions of smooth threefolds are divisorial the exceptional locus
of $\mu_i$ is exactly $D_i$.
Since $D_i \cdot l = 0$, a general
deformation of $l$ does not meet the exceptional 
locus, so the morphisms $\mu_i$ are equivalences
of models.

By the classification of birational Mori extremal contractions 
in dimension $3$, the morphism $\mu_i$
contracts $D_i$ to a point, or $X$ is the blow-up of 
a manifold $X_i$ along a smooth curve.
The model $X'$ is obtained along the following algorithm.

{\it Step 1.} One of the morphisms $\mu_i$ is the blow-up 
of a Fano manifold $X_i$ along a smooth curve.
We replace $(X,l)$ by the equivalent model $(X_i,\mu_i(l))$ and restart the program.

{\it Step 2.} None of the morphisms $\mu_i$ is the 
blow-up of a Fano manifold $X_i$ along a smooth curve.
It follows from \cite[Prop.4.5.]{MM83} 
and \cite{Mo82} that for all $i$ the line bundle
$\sO_{D_i}(D_i)$ is negative on all the curves contained in $D_i$.
More precisely, either $D_i$ is contracted to a point,
or $D_i$ is isomorphic to $\PP^1\times \PP^1$
with normal bundle 
$\sO_{D_i}(D_i)\simeq {\mathcal O}_{\PP^1 \times \PP^1}(-1,-1)$
and $D_i$ is mapped to a smooth rational curve with normal bundle
${\mathcal O}_{\PP^1}(-1)\oplus {\mathcal O}_{\PP^1}(-1)$
(and in that case $X_i$ is smooth). 
A standard argument shows that $D_i \cap D_j= \emptyset$ for all $i \neq j$.
We can then ``do successively'' the morphisms 
$\mu_i$, getting a birational morphism 
\holom{\mu}{X}{X'} such that the exceptional divisor
is exactly $D_1 \cup \ldots \cup D_k$. Moreover the variety $X'$ is a 
$\Q$-Gorenstein projective threefold 
with at most terminal singularities since at each step, the morphism
$\mu_i$ is a Mori contraction.

Finally, 
since a general quasi-line $l$ does not meet $D_1 \cup \ldots \cup D_k$, 
the morphism $\mu$ is an equivalence 
of models $(X,l) \simeq (X',\mu(l))$. Since the strict transform of 
any effective prime divisor on $X'$ is an effective prime divisor on $X$, 
it is 
also clear
that $X'$ is minimal with respect to $\mu(l)$. \finpreuve

\section{Fibrations between models}
\label{sectionfibrations}

Let $(X,l)$ be a smooth model  
that admits a fibration $\holom{\varphi}{X}{Y}$. 
Since $l$ is a very free curve,
it is clear that $l$ is not contracted by $\varphi$
and the image $\varphi(l)$ is also a very free curve in $Y$. 
The answer to the following questions
is not so clear: 
\begin{enumerate}
\item Do $Y$ or the general fibre $F$ contain quasi-lines ? How are they related to $l$ ? 
\item Does the minimality with respect to $l$ imply any restrictions on the fibre space structure?
\end{enumerate}
In general one should not expect too much, but we will see that the  restriction to
morphisms of models $(X,l) \rightarrow (Y,\varphi(l))$ is a good framework for analyzing these problems.

\begin{lemm}
\label{lemmafibrations}
Let \holom{\varphi}{X}{Y} a fibration from 
a projective manifold $X$ onto a projective manifold $Y$ of dimension 
at least 2. 
\begin{enumerate}
\item Let $l$ be
a general deformation of a given quasi-line on 
$X$ and suppose that $\varphi(l)$ is smooth.
Then the following are equivalent:
\begin{enumerate}
\item the image $\varphi(l) \subset Y$ is a quasi-line,
\item the curve 
$l \subset \fibre{\varphi}{\varphi(l)}$ is a quasi-line
of the manifold $\fibre{\varphi}{\varphi(l)}$. 
\end{enumerate} 
In case (a) or (b), the morphism $l \rightarrow \varphi(l)$ is an isomorphism.
\item Conversely, let $l$ be a rational curve in $X$ such that 
the image $\varphi(l) \subset Y$ is a quasi-line,
the rational curve $l \subset \fibre{\varphi}{\varphi(l)}$ is a quasi-line
of the manifold $\fibre{\varphi}{\varphi(l)}$ and 
the morphism $l \rightarrow \varphi(l)$ is an isomorphism. Then
$l$ is a quasi-line of $X$.
\end{enumerate}
\end{lemm}

{\em Proof.} 
Note first that since $\varphi(l)$ is general and smooth, 
the preimage $Z:=\fibre{\varphi}{\varphi(l)}$ is smooth 
of dimension $\dim X - \dim Y + 1$.
So we have an exact sequence
\[
0 \rightarrow N_{l/Z} \rightarrow N_{l/X} 
\rightarrow N_{Z/X} \otimes \sO_l \rightarrow 0,
\]
and $N_{Z/X} \simeq \varphi^* N_{\varphi(l)/Y}$.

The second statement is then obvious since any extension 
\[ 
0 \rightarrow {\mathcal O}_{\PP^1}(1) ^{\oplus a}   
\rightarrow E \rightarrow {\mathcal O}_{\PP^1}(1) ^{\oplus b} \rightarrow 0
\]
splits.

For the first statement, the isomorphism 
$N_{Z/X} \simeq \varphi^* N_{\varphi(l)/Y}$ implies
\[
\deg N_{Z/X} \otimes \sO_l = d \cdot \deg  N_{\varphi(l)/Y},
\]
where $d$ is the degree of the morphism $l \rightarrow \varphi(l)$.
We have $\deg N_{l/X}= \dim X - 1$ and $\varphi(l)$ is a very free curve, so 
$\deg N_{\varphi(l)/Y} \geq \dim Y - 1 \geq 1$. Hence
\begin{eqnarray*}
(*) \qquad \deg N_{l/Z} &=& \deg N_{l/X} - \deg_l N_{Z/X} \otimes \sO_l
\\
&=& \dim X - 1 - d \cdot \deg  N_{\varphi(l)/Y}
\\
&\leq& \dim X - 1 - d \cdot (\dim Y - 1)
\\
& \leq& \dim X - 1 - \dim Y + 1 = \dim Z - 1,
\end{eqnarray*}
and equality holds if and only if $d=1$ and 
$\deg N_{\varphi(l)/Y} = \dim Y - 1$. 

$(b) \Rightarrow (a)$. If $l \subset Z$ is a quasi-line, 
then $\deg N_{l/Z}  = \dim Z-1$, so equality holds in $(*)$.
Hence $\deg N_{\varphi(l)/Y} = \dim Y - 1$ and $d=1$, so
$\varphi(l)$ is a quasi-line and the
morphism  $l \rightarrow \varphi(l)$ is an isomorphism.

$(a) \Rightarrow (b)$.
Suppose now that the image of a general quasi-line of $X$ 
is a quasi-line of $Y$. 
Fix a general point $x$ in $X$ such that for some quasi-line $l$
through $x$, the image $\varphi(l)$ is a quasi-line through $\varphi(x)$.
Let $\mathcal H_x$ and 
$\mathcal H'_{\varphi(x)}$ the corresponding varieties as defined in 
section \ref{sectionbasics}. Since $\varphi(l)$ is a quasi-line we have by \cite[I,Thm.6.8]{Ko96} a natural surjective morphism
\holom{\overline{\varphi}}{\mathcal H_x}{\mathcal H'_{\varphi(x)}} defined by $[l] \mapsto [\varphi_*(l)]$.
Choose $[l']$ such that
$\fibre{\overline{\varphi}}{[l']}$ has an irreducible component
$S$ of dimension $\dim X - \dim Y$ 
and 
parametrises at least one quasi-line. We may also assume that 
$Z:=\fibre{\varphi}{l'}$
is a smooth variety. By construction,
$p(\fibre{q}{S}) \subset Z$ 
(we use the notation of diagram \ref{basicdiagram}).
By construction $S$ parametrises a family of cycles through a fixed point $x$ that dominates $Z$,
since a general member $l \in S$ is irreducible, it is a very free curve in $Z$.
This shows that $\deg N_{l/Z} \geq \dim Z - 1$. 
By the inequality $(*)$ we even have equality, 
so $l$ is a quasi-line in $Z$. $\finpreuve$

\medskip

{\bf Remark.} The proposition shows that 
the study of fibrations \holom{\varphi}{(X,l)}{(Y,\varphi(l))},
where $l$ and $\varphi(l)$ are quasi-lines, 
essentially reduces to the study of fibrations
\holom{\varphi}{(X,l)}{\PP^1} where $l$ is both a quasi-line and a 
section of $\varphi$. Note however that in general, 
$\fibre{\varphi}{\varphi (l)}$ is not Fano even if $X$ and $Y$ are (see 
example (2) immediately after theorem \ref{theoremnotconicbundle}).
One question also immediately arises:
if $l$ is a section of \holom{\varphi}{(X,l)}{\PP^1}, 
does the general fibre contain a quasi-line?
We will show in the next paragraph
that this is true when $X$ is a Fano manifold
of dimension at most
four, but the problem becomes more delicate
in higher dimension.

\subsection{Fibrations over a curve}
\label{subsectionfibrationcurve}

Recall that the pseudo-index 
of a Fano manifold $X$ is the positive integer defined 
as $$i_X := \min 
\{ -K_X \cdot C \,|\, C \mbox{ is a rational curve of } X\}.$$
This number has been very much studied in the last period, the 
general philosophy being that the Fano manifolds with high pseudo-index
are the easiest to understand.

\begin{lemm}
\label{lemmapseudoindex}
Let $(X,l)$ be a smooth Fano model of dimension $n$,  
and let $\varphi : X \to \PP^1$ be a fibration such that $l$ is a section of $\varphi$. 
Then one of the following holds:
\begin{enumerate}
\item the pseudo-index of the general fibre is strictly smaller than $(n+1)/2$, 
\item the general fibre of $\varphi$ is isomorphic to $\PP^{n-1}$.
\end{enumerate}
If moreover $X$ is minimal with respect to $l$, the second case does not occur.
\end{lemm}

{\em Proof.} Fix a general point $x \in X$ such that 
$F:= \varphi ^{-1}(\varphi(x))$ is a smooth fibre and such that there passes
a quasi-line through $x$. Let $y \in F$ be  a point different from $x$,
then there exists no quasi-line through $x$ and $y$: the quasi-line 
is a section, so it meets $F$ in exactly one point.

Therefore there 
exists a degeneration $\sum_{i=1}^{k} \alpha_i l_i$ connecting $x$ to $y$.
Then we have $F \cdot l_i  \geq 0$ for all $i$,  
and equality holds if and only if $l_i$ is contracted by $\varphi$.
Since $l$ is a section, we have
\[
1 = F \cdot l = F \cdot \sum_{i=1}^{k} 
\alpha_i l_i = \sum_{i=1}^{k} \alpha_i F \cdot l_i.
\]
It follows that there exists a unique component, 
say $l_1$ such that $F \cdot l_1 = 1$ and for $i>1$, we have
$F \cdot l_i = 0$.
Since for $i>1$, the curves $l_i$ are contained in a $\varphi$-fibre,
there exists a renumbering of $l_2, \ldots l_k$ and  
$2 \leq k' \leq k$ such that for
$2 \leq i \leq k'$, the curve $l_i$ is contained in $F$ 
and the connected chain of curves $l_2 \cup \ldots \cup l_{k'}$ passes through $x$ and $y$.

We distinguish two cases: 
\begin{enumerate}
\item[a)] there exists a $y \in F$, such that $k'>2$. Denote by 
$i_F$ the pseudo-index of $F$, then
\[
n+1 = - K_X \cdot l = - K_X \cdot (l_1 + \sum_{i=2}^{k} \alpha_i l_i) 
> - K_X \cdot \sum_{i=2}^{k'} \alpha_i l_i \geq 2 i_F, 
\]
so we are in the first case of our statement,

\item[b)] for all $y \in F \setminus \{x\}$, we have $k' = 2$. 
Then the curves $l_2$ form an unsplit family of curves that
connects every point to $x$. 
It follows, again by theorem \ref{theoprojectivespace} that $F \simeq \PP^{n-1}$.
Assume now that $X$ is minimal with respect to $l$. Then every 
fibre of $\varphi$ is irreducible.
Since $l$ is a section, every 
fibre of $\varphi$ is reduced. Hence 
$X$ is a projective bundle over $\PP^1$
(by a simple refinement, due to Araujo, of Fujita's characterisation
of projective bundles \cite{Fuj87}). 
Finally  
we have $X  \simeq \PP(\sO _{\mathbb P _1} ^{\oplus n-1} \oplus \sO  _{\mathbb P _1}(1))$
by \cite[Prop.4.1]{IV03}, 
which is not minimal.
\end{enumerate}
~\finpreuve

\medskip

The next example shows that the lemma does not hold if the quasi-line is not a section.

\medskip

{\bf Example.} 
Let $X$ be a general member of 
${\mathcal O}_{\PP^1 \times \PP^n}(1,d)$, where 
$n,d$ are
integers satisfying $n \geq 3$ and $1 \leq d \leq n $.  
Then $X$ is Fano and the general fibre of $\varphi: X \to \PP^1$
is a smooth hypersurface of $\PP^n$ 
of degree $d$ (hence with pseudo-index $n+1-d$).
Moreover, $X$ contains a quasi-line $l$ (the inverse image of a general line
in $\PP^n$) but the map $\varphi_{|l}: l\to \PP^1$ has degree $n$. 

\begin{theo}
Let $(X,l)$ be a smooth Fano model 
of dimension at most four, and let $\varphi: X \to \PP^1$ be a fibration 
such that $l$ is a section of $\varphi$.
Then the general fibre $F$ contains a quasi-line.
\end{theo}

{\em Proof.} The statement is trivial for $\dim X=1,2$.
If $\dim X=3$, the general fibre
is a del Pezzo surface. All the del Pezzo surfaces except $\PP^1 \times \PP^1$
contain quasi-lines, but the quadric is excluded by lemma \ref{lemmapseudoindex}.
If $\dim X=4$, the result is an immediate consequence of the lemma
\ref{lemmafamily} below which gives a much more precise information in this case:
the quasi-line of $F$ is a component of a degeneration of $l$. \finpreuve 

\medskip

Before stating a technical lemma, recall that for any smooth 
quasi-projective variety $Y$, there exists a subset $Y^{\rm free}$
which is the intersection of countably many dense open subsets of $Y$ 
such that any rational curve on $Y$ whose image meets $Y^{\rm free}$
is free.   

\begin{lemm}
\label{lemmaFanochain}
Let $F$ be Fano manifold of dimension at most three
and let $x \in F$ be
a fixed very general point.
Assume that there is a connected chain
of rational curves $l'_1 \cup \ldots \cup l'_k$ passing through $x$
whose deformations cover a dense open subset of $F$ and such
that $-K_F \cdot \sum_{i=1}^k l'_i \leq  \dim F+1$.
Then a general deformation satisfies $k=1$ and $l'_1$ is a very free rational curve in $F$ such that 
$-K_F \cdot l'_1=\dim F+1$.
\end{lemm}

{\bf Remark.} This lemma is the key ingredient for lemma \ref{lemmafamily}, and it is at this point
that we use the hypothesis on the dimension. The authors of this paper are not convinced that the statement
generalises to arbitrary dimension.

\medskip

{\em Proof.} The statement is trivial is $F$ is a curve, so suppose $\dim F$ is $2$ or $3$.
We argue by contradiction and assume that for 
$y \in F$ a very general point,
the connected chain of curves $\cup_{i=1}^{k} l'_i$ passing through $x$ and $y$ is reducible,
in particular $k \geq 2$. Note that since $x$ and $y$ are very general they are both in $F^{\rm free}$.
Furthermore we can suppose that no irreducible component $l'_{i_0}$
contains both $x$ and $y$; otherwise the corresponding component $l'_{i_0}$ would be very free
\cite[Prop.4.20]{Deb01}, so
\[
\dim F+1 = -K_F \cdot \sum_i l'_i > -K_F \cdot l'_{i_0} \geq \dim F+1,
\]
a contradiction. 

So up to renumbering we can suppose that $x \in l'_1$, $y \in l'_2$. 
Since $x \in F^{\rm free}$, we have $-K_F \cdot l'_1 \geq 2$.
Since $y \in F^{\rm free}$, we also have $-K_F \cdot l'_2 \geq 2$.
If $\dim F=2$, this contradicts $3 \geq -K_F \cdot \sum_i l'_i$.
Suppose now that $\dim F=3$, then 
\[
4 \geq - K_F \cdot \sum_{i=1}^{k} l'_i,
\]
so the preceeding inequalities imply $k=2$ and
$-K_F \cdot l'_1 = - K_F \cdot l'_2= 2$. 
Since $x \in F^{\rm free}$ there exist only finitely many rational curves $C$
such that $-K_F \cdot C = 2$ and $x \in C$. 
Therefore there are at most finitely many
curves $l'_1$ passing through $x$, we denote their union by $L'_1$.
Since $l'_2$ is free and $-K_F \cdot l'_2 = 2$, 
we have $N_{l'_2/F} \simeq \sO_{\PP^1}^{\oplus 2}$.
It follows that $\Hilb(F)$ is smooth of dimension two at the point $[l'_2]$
and we choose a smooth open neighbourhood $Z \subset \Hilb(F)$
parametrizing deformations of $l'_2$ with the same normal bundle.
Let $\Gamma$ be the universal family over $Z$, and 
let \holom{p}{\Gamma}{F} and \holom{q}{\Gamma}{Z} be the natural maps. 
Then $p$ is finite over its image 
and since $\dim L'_1=1$ and $\dim Z=2$, this implies
$q(\fibre{p}{L'_1} \subsetneq Z$. Therefore a general chain of curves $l'_1 \cup l'_2$ is not connected, 
a contradiction.
\finpreuve

\begin{lemm}
\label{lemmafamily}
In the situation of the theorem, let $x \in X$ be 
a fixed very general point such that $x \in F^{\rm free}$ for $F=\fibre{\varphi}{\varphi(x)}$ a smooth
fibre. Then there exists
a quasi-projective $(\dim X-2)$-dimensional family of
degenerations of $l$ with fixed point $x$ that are of the form
\[
l_1 + l_2,
\]
where $l_1$ is a section of $\varphi$ such that 
$-K_X \cdot l_1 =1$ and $l_2 \subset F = \varphi ^{-1}(\varphi(x))$
is a very free rational curve in $F$ such that $-K_F \cdot l_2=\dim F+1$.
If $\dim F=3$, the curve $l_2$ is a quasi-line of 
$F$\footnote{By \cite[II,Thm.3.14]{Ko96} a general deformation
of a very free curve is an embedding if the ambient variety has dimension at least three, but only an immersion
if the dimension is two. Recall that we defined quasi-lines to be smooth curves.}. 
\end{lemm}

{\em Proof.} We assume that $\dim X =3$ or $4$, the other cases being trivial.

\smallskip

As seen in the proof of lemma \ref{lemmapseudoindex}, 
for all $y \in F \setminus \{x\}$ 
the degenerations of $l$ through $x$ 
that connect $x$ to $y$
are of the form
\[
l_1 + \sum_{i=2}^{k'} \alpha_i l_i +\sum_{i=k'+1}^{k} \alpha_i l_i ,
\]
where $l_1$ is a section of $\varphi$, the connected chain of curves 
$\sum_{i=2}^{k'} \alpha_i l_i$ is contained in $F$ and 
passes through $x$ and $y$.
By lemma \ref{lemmaFanochain} 
the connected chain of curves
$\cup_{i=2}^{k'} l_i$ passing through $x$ and $y$ is irreducible,
so $k'=2$ and $l_2$ contains $x$ and $y$.
By the same lemma 
for $y \in F$ very general, 
the curve $l_2$ is very free in $F$ and $-K_F \cdot l_2 = \dim F+1$.
Since $y$ varies in a family whose closure has dimension $\dim X-1$, we obtain a 
$(\dim X-2)$-dimensional family of degenerations having 
the stated structure. \finpreuve

\begin{theo}
Let $(X,l)$ be a smooth Fano model 
of dimension at most four, and let $\varphi: X \to \PP^1$ be a fibration 
such that $l$ is a section of $\varphi$.
If $e(X,l)=1$, then the variety $X$ is not minimal with respect to $l$.
\end{theo}

{\em Proof.} The case $\dim X=2$ is treated in example \ref{examplesurfaces},
suppose now that $\dim X \geq 3$.
Fix a general point $x \in X$ and consider the closure in $\chow{X}$
of the family given by lemma \ref{lemmafamily}. We obtain a  
$(\dim X-2)$-dimensional family $S$ of
degenerations of $l$ 
with fixed point $x$ such that the general degeneration is of the form
$l_1 + l_2$,
where $l_1$ is a section of $\varphi$ with 
$-K_X \cdot l_1 =1$  and $l_2 \subset F = \varphi ^{-1}(\varphi(x))$
is a very free rational curve in $F$ with $-K_F \cdot l_2=\dim F+1$.
Let $\Gamma \subset S \times X$ be the normalization of the incidence variety, and 
let \holom{p}{\Gamma}{X} and \holom{q}{\Gamma}{S} be the natural maps. 
By construction 
the fibre $F:=\fibre{\varphi}{\varphi(x)}$ is contained in $p(\Gamma)$, 
but does not contain the curves $l_1$. 
Since $F$ is a divisor this shows that $p(\Gamma)$ 
has at least two irreducible components,
that is the fibre $F$ and a component $D$ that contains the curves $l_1$.
We claim that $D$ is a divisor that is dominated by the deformations of $l_1$.
Assume this for moment, then 
$-K_X \cdot l_1=1$ implies that the deformations of $l_1$ 
form an unsplit family, so they cover $D$. This implies that
$x \notin D$, otherwise there would exist a non-free curve passing through $x$.   
Conclude with lemma \ref{lemmabirationalevaluation}.

{\em Proof of the claim.}
We argue by contradiction and suppose that the locus $V \subset D$ dominated by the curves
$l_1$ is not a divisor.
Since $V \not\subset F$, the intersection $V \cap F \subset F$ is not a divisor. 
Since the chain  
$l_1 \cup l_2$ is
connected and $l_1 \cap l_2 \subset V \cap F$, all the curves 
$l_2$ pass through
$V \cap F$ and $x$. This is impossible by \cite[II, Prop.3.7]{Ko96}. 
\finpreuve

\medskip

{\bf Remark.} The proof above uses the hypothesis on the dimension of $X$ 
only to assure the existence of a family of degenerations of the form $l_1+l_2$.
\medskip

The next proposition shows how one can recover a quasi-line by smoothing 
degenerations of the form $l_1+l_2$.

\medskip

\begin{prop}
\label{propositionsmoothing}
Let $X$ be a Fano manifold of dimension $n$, and let $\varphi: X \to \PP^1$ be an elementary contraction 
such that the general fibre $F$ contains a quasi-line $l_2$, that is 
\[
T_F|_{l_2} \simeq \sO_{\PP^1}(2) \oplus \sO_{\PP^1}(1)^{\oplus n-2}.
\]
Suppose furthermore that the second elementary contraction $\psi~: X \to Y$
has length 1, 
and that a generator $l_1$ of the corresponding extremal ray 
is a section of $\varphi$.
Then $X$ contains a quasi-line that is a section of $\varphi$.
\end{prop}

{\em Proof.} 

{\em Step~1.}
Since the second contraction $\psi$ has length 1, it is either a conic bundle 
(with singular fibre) or 
$X$ is a blow-up of the smooth variety $Y$ along a smooth subvariety of codimension two \cite[Cor.1.4]{Wi91}. In both cases, 
there exists an irreducible divisor $D \subset X$ that is dominated by the deformations of
a generator $l_1$ of the second extremal ray and
\[
T_X|_{l_1} \simeq \sO_{\PP^1}(2) \oplus \sO_{\PP^1}^{\oplus n-2} \oplus \sO_{\PP^1}(-1).
\]
Fix now a general point $x \in X \setminus D$ such that 
$F=\fibre{\varphi}{\varphi(x)}$ is smooth and the quasi-lines
$l_2$ passing through $x$ dominate $F$. 
Since $\varphi$ is an elementary contraction 
of fibre type, the divisor $D$ is strictly positive
on the first extremal ray, in particular $(D \cap F) \cdot l_2=D \cdot l_2>0$.
Hence for a general curve $l_1$ there exists a quasi-line 
$l_2$ such that $C:=l_1 \cup l_2$ is connected. 
Since $l_1$ and $l_2$ are smooth curves, the
intersection $l_1 \cap l_2 = l_1 \cap F$ is a reduced point.
The connected curve $C$ is a tree of rational curves and
we will follow closely the argument in \cite[Prop.4.24]{Deb01}
to show that $C$ is smoothable fixing the point $x$ 
(we refer to \cite[Ch. 2]{Deb01} for the notation
and details on schemes parametrizing morphisms).
Consider the inclusion morphism 
\holom{f}{C}{X} and let $b \in l_2$ be the unique point
such that $f(b)=x$.
Let \holom{\pi}{\mathcal C}{(T,0)} be a 
smoothing\footnote{Such a smoothing always exists and we can take $T$ to 
be the unit disc in $\C$, cf. \cite[p.101]{Deb01}.} of the tree $C$ 
and let \holom{\sigma}{T}{\mathcal C}
be a section such that $\sigma(0)=b$. 
Denote by \holom{g}{\sigma(T)}{X \times T}
the morphism defined by $g(\sigma(t))=(b,t)$.
The $T$-morphisms from $\mathcal C$ to 
$X \times T$ extending $g$ are parametrised by
the $T$-scheme ${\rm Mor}_T(\mathcal C,X \times T; g)$ 
whose fibre at $0$ is ${\rm Mor}(C,X,b\mapsto x)$.
By \cite[II,Thm.1.7]{Ko96} the irreducible components of 
${\rm Mor}_T(\mathcal C, X \times T; g)$ 
at $[f]$ have dimension at least
\[
\chi(C, f^* T_{X}(-b)) + \dim T.
\]

\smallskip

{\em Step 2.} 
We claim that ${\rm Mor}(C,X,b\mapsto x)$ has a component of the expected dimension, that is equal to
$\chi(C, f^* T_{X}(-b)).$ 

\smallskip

The exact sequence
\[
0 \rightarrow f^* T_X(-b-(l_1 \cap l_2)) \otimes \sO_{l_2} 
\rightarrow f^* T_X(-b)
\rightarrow f^* T_X(-b) \otimes \sO_{l_1}
\rightarrow 0.
\]
immediately implies $\chi(C, f^* T_{X}(-b)) = n+1$.
Moreover
the general fibres of the restriction morphism
${\rm Mor}(C,X,b\mapsto x) \to {\rm Mor}(l_1,X)$
consist of quasi-lines of $F$ passing through
two fixed points, hence are $0$-dimensional.
Therefore there exists at least one component of $Z \subset {\rm Mor}(C,X,b\mapsto x)$ 
that dominates ${\rm Mor}(l_1,X)$ and its dimension is equal to
$$ \dim {\rm Mor}(l_1,X) = h^0(l_1, f^* T_X \otimes \sO_{l_1})
= h^0(l_1,\sO_{\PP^1}(2) 
\oplus \sO_{\PP^1}^{\oplus n-2} \oplus \sO_{\PP^1}(-1))=n+1.$$

\smallskip

Therefore the fibre at $0$ of the morphism ${\rm Mor}_T(\mathcal C,X \times T; g) \rightarrow T$ 
has a component of dimension $\chi(C, f^* T_{X}(-b))$.
Since the components of ${\rm Mor}_T(\mathcal C,X \times T; g)$  have dimension at least
$\chi(C, f^* T_{X}(-b)) + \dim T$,
it follows that ${\rm Mor}_T(\mathcal C,X \times T; g)$ has a component of 
the expected dimension $\chi(C, f^* T_{X}(-b)) + \dim T$ that dominates $T$, hence
the morphism ${\rm Mor}_T(\mathcal C,X \times T; g) \rightarrow T$ 
is dominant at $[f]$.
In particular there exists in 
${\rm Mor}_T(\mathcal C,X \times T; g)$ an irreducible (open) curve 
passing through $[f]$ that dominates $T$. Denote by $T'$ its normalization, then the morphism
$T' \rightarrow {\rm Mor}_T(\mathcal C,X \times T; g)$  
yields a $T'$-morphism 
\[
\mathcal C \times_T T' \rightarrow X \times T'
\]
which is the smoothing to a rational curve $l$ keeping fixed $f(b)=x$.

\smallskip

{\em Step~3.} 
Clearly $-K_X \cdot l= n+1$ and $F \cdot l= F \cdot l_1 + F \cdot l_2=1$.
If we show that the deformations of $l$ keeping fixed $x$ dominate $X$, we
conclude by \cite[Prop. 4.20]{Deb01} that a general deformation of $l$ 
is very free, so it is a quasi-line.
We argue by contradiction and suppose that this is not the case. 
Then the locus $V$ of deformations and degenerations of $l$ 
is a finite union
of proper subvarieties, and the divisors 
$F$ and $D$ are irreducible components of this locus. Choose now 
a connected cycle $l_1 \cup l_2$ that is not contained in any irreducible component of $V$. 
Since $l_1 \cup l_2$ is smoothable by some deformation $l$, 
we can choose an irreducible curve $l'$ arbitrarily close to 
$l_1 \cup l_2$ that is contained in $V$.
Yet this implies that $l' \subset F \cup D$, so by irreducibility 
of $l'$ we have
$l' \subset F$ or $l' \subset D$. Since neither $l_1 \cup l_2 \subset F$ 
nor  $l_1 \cup l_2 \subset D$
we obtain a contradiction. \finpreuve

\subsection{New examples} \label{subsectionexamples}
{\em In this section only,} $X$ is a Fano threefold
with Picard number $2$ which has a fibration
$\varphi~: X \to \PP^1$ and a quasi-line $l$ being a section of $\varphi$.
Since $X$ has Picard number $2$, there is another
Mori extremal contraction $\psi~: X \to Y$, whose fibres have dimension
at most one, hence either $\psi$ is a conic bundle, or $\psi$
is the blow-up of the smooth variety $Y$ with centre
a smooth curve. We will deal with the case of conic bundles in the subsection
\ref{subsectionconicbundles}. 
In the case where $\psi$
is the blow-up of a smooth variety $Y$ with centre
a smooth curve, the situation is very nice: $\psi$ is determined
by the degenerations of $l$.

\begin{cor}
Let $(X,l)$ be a smooth Fano model of dimension three that admits an 
elementary contraction $\varphi: X \to \PP^1$  such 
that $l$ is a section of $\varphi$.
Assume that the second elementary Mori contraction is a birational contraction $\psi: X \to Y$. 
Then both extremal rays of ${\rm NE}(X)$ are spanned
by irreducible component of degenerations of $l$.
\end{cor}

{\em Proof.} By lemma \ref{lemmafamily}
there exists a positive-dimensional quasi-projective 
family of degenerations of $l$ 
with fixed point a general point $x$ that are of the form
\[
l_1 + l_2,
\]
where $l_1$ is a section of $\varphi$ such that 
$-K_X \cdot l_1 =1$ and $l_2 \subset F = \varphi ^{-1}(\varphi(x))$
satisfies $-K_F \cdot l_2=3$. In particular, the extremal ray 
defining $\varphi$ is spanned by $l_2$.  

Let $C$ be a generator of the extremal ray contracted by $\psi$, 
then $F \cdot C=1$ by \cite[Prop.6]{MM81}. 
Let $E$ be the exceptional divisor of $\psi$, then 
there exists rational numbers $a$ and $b$ such that
$E = a (-K_X) + b F$. Hence

\begin{eqnarray*}
E \cdot C = -1 &=& a + b
\\
E \cdot l &=& 4 a  + b
\\
E \cdot l_2 &=&  3 a. 
\end{eqnarray*}

This implies 
$E \cdot l_1 = E \cdot l - E \cdot l_2 = a+b= -1 = E \cdot C $.
Since furthermore 
$ -K_X  \cdot l_1 = 1 = -K_X  \cdot C$ and $X$ has Picard number $2$,
this implies that $l_1$ and $C$ have the same numerical class. 
It follows that $l_1$ is contracted by $\psi$,
so it generates the second extremal ray.\finpreuve

\medskip

{\bf Examples.} Under the assumptions of the corollary,
and using again \cite{MM81,MM83}, the threefold $Y$ is Fano of index $r \geq 2$
and easy computations show that $a = r-1$, therefore $b= -r$
and $E \cdot l = 3r-4$. Moreover, since $E \cdot l_1 = 1$,
$\psi |_F~: F \to F'=\psi(F)$ is an isomorphism for any
general $\varphi$-fibre $F$ and $F' \in |(r-1){\mathcal O}_{Y}(1)|$.  
Moreover, we know (by \cite{MM81,MM83} again) that $Y$ is either
$\PP^3$, ${\mathcal Q}_3$ or one of 
the five Fano threefolds $V_d$ ($1\leq d \leq 5$)
of index $2$ and Picard number $1$. Here are three explicit examples:
\begin{enumerate} 
\item[(1)] Let $\psi~: X \to \PP^3$ be the blow-up of $\PP^3$ with centre
the transverse intersection $C = S_1 \cap S_2$ of 
two smooth cubic surfaces $S_1$ and $S_2$. 
Then the strict transform $l$ of a normal twisted rational 
cubic curve meeting $C$ transversaly in exactly $8$ points
is a quasi-line of $X$, and a section of the projection
$\varphi~: X \to \PP^1$ corresponding to the pencil defined by $S_1$ and $S_2$. 
 
\item[(2)] Let $\psi~: X \to {\mathcal Q}_3 \subset \PP^4$ 
be the blow-up of the $3$-dimensional quadric ${\mathcal Q}_3$ with centre
the transverse intersection $C = S_1 \cap S_2$ of 
two smooth members $S_1$ and $S_2$ of $|{\mathcal O}_{{\mathcal Q}_3}(2)|$. 
Then the strict transform $l$ of a normal twisted rational 
cubic curve meeting $C$ transversally in exactly $5$ points  
is a quasi-line of $X$, and a section of the projection
$\varphi~: X \to \PP^1$ corresponding to the pencil defined by $S_1$ and $S_2$.

\item[(3)] Let $V_3$ be a smooth cubic hypersurface of $\PP^4$, and let
$\omega \subset V_3$ be a normal twisted rational cubic contained
in $V_3$. Denote by $P_3$ the hyperplane in $\PP^4$ generated by $\omega$. 
Let finally $l'$ be a chord of $\omega$ ({\em i.e.} a 
line in $\PP^4$ meeting $\omega$ in two points)
and $P_1$ and $P_2$ be two general hyperplanes  
in $\PP^4$ such that $l'= P_1 \cap P_2 \cap P_3$. Let $\psi~: X \to V_3$
be the blow-up of $V_3$ with centre $C = V_3 \cap P_1 \cap P_2$ ($C$
is therefore a plane cubic).
Since $\omega \cap C \subset V_3 \cap l' \subset \omega \cap l'$, 
the curve $\omega$ meets $C$ in exactly two points. 
The strict transform of $\omega$   
is a quasi-line of $X$, section of the projection
$\varphi~: X \to \PP^1$ corresponding to the pencil defined by 
$V_3 \cap P_1$ and $V_3 \cap P_2$.

\end{enumerate} 

\subsection{Conic bundles}
\label{subsectionconicbundles}

The variety $\PP(T_{\PP^2})$ is a Fano threefold
that contains a quasi-line $l$ such that the natural projection to 
$\PP^2$ induces a morphism 
of models $(\PP(T_{\PP^2}),l) \rightarrow (\PP^2,line)$. 
One can obtain conic bundles with the same property by blowing-up 
$\PP(T_{\PP^2})$  along a smooth curve, 
but such a manifold would of course not be minimal. We will 
now show that this is what always happens.

\medskip
 
\begin{theo}
\label{theoremnotconicbundle}
Let \holom{\varphi}{(X,l)}{(Y,\varphi(l))} be a morphism of smooth models,
and suppose that $\varphi$ is flat of relative dimension $1$.
Assume furthermore that:
\begin{enumerate}
\item[(1)] $X$ is minimal with respect to $l$ and $e(X,l)=1$,
\item[(2)] there exists a big and nef line bundle $H$ on $Y$ 
such that $\varphi^* H \cdot l=1$. 
\end{enumerate}
Then $Y\simeq \PP^{n-1}$ and $X \simeq \PP(E)$, where $E$ is a rank two vector bundle over $Y$.
\end{theo}

{\em Proof.} Since $X$ is minimal with respect to $l$, the manifold 
$Y$ is minimal with respect to $\varphi (l)$. 
By assumption (2) and 
theorem \ref{theoremalmostline}, 
we get that $Y \simeq \PP^{n-1}$ and $\varphi(l)$ 
is a line. 
Suppose first that $\varphi$ is a smooth morphism, so all the fibres are isomorphic to $\PP^1$.
Therefore $\varphi$ is locally analytically trivial and since the Brauer group of $\PP^{n-1}$ vanishes,
we see that $X$ is the projectivisation of a rank two bundle.  

\smallskip

We will now argue by contradiction and suppose that $\varphi$ is not a smooth morphism. 
We claim that in this case there exists 
an effective divisor $\Delta$ in $Y$ such that the general
$\varphi$-fibre over $\Delta$ is reducible. 
Since the general fibre is a rational curve \cite[Thm.1.12]{IV03}, this is clear if 
$\varphi$ is elementary (that is $\rho(X) - \rho (Y) =1$): $\varphi$
is a conic bundle and the discriminant locus is such a $\Delta$. 
If $\rho(X) - \rho (Y) > 1$, just apply the relative contraction theorem 
to get  an elementary contraction \holom{\mu}{X}{X'}
that is a $Y$-morphism. It is not hard to see that $\mu$ 
is birational and by flatness of $\varphi$, 
all the fibres of $\mu$ have dimension at most $1$,
so $X'$ is smooth and, by Ando's theorem \cite[Thm.2.1]{An85}, 
the morphism $\mu$ is the blow-up of a smooth codimension $2$ 
subvariety $Z$ such that
the image of $Z$ in $Y$ is a divisor. In particular there 
exists a divisor $\Delta \subset Y$
such that the general fibre over $\Delta$ is reducible.

\smallskip

Fix now an irreducible divisor $\Delta \subset Y$ 
such that the general $\varphi$-fibre over $\Delta$ is reducible,
fix also a general point $x \in X$ such that $x \notin \fibre{\varphi}{\Delta}$. 
We consider the deformations of $l$ with fixed point $x$. Recall the notation
of the basic diagram \ref{basicdiagram}. 
\[
 \xymatrix{ 
{\mathcal U}_x \ar[d]^{q}\ar[r]^{p} & X \\
{\mathcal H}_x & } 
\]
By lemma \ref{lemmabirationalevaluation}, 
we have reached a contradiction if we show that an irreducible 
component $D \subset \fibre{\varphi}{\Delta}$ 
is covered by irreducible components of degenerations with fixed point $x$.

\medskip

Let $y$ be a general point of $\Delta$, 
then the fibre $\fibre{\varphi}{y}$ 
has at least one singular point $z$ at the intersection of two components.
Since $\varphi^* H \cdot l=1$, every quasi-line $l$ parametrised 
by $\mathcal H_x$ is a subsection of $\varphi$, 
in particular it does not meet the singular points of fibres. 
Since the map $p$ is surjective,
this shows that there exists a degeneration 
$l \sim \sum_{i =1}^k \alpha_i l_i$ that connects $x$ to $z$.
Since $\varphi^* H \cdot \sum_{i =1}^k \alpha_i l_i =1$, this cycle 
can be written
$l_1 + \sum_{i \geq 2}^k \alpha_i l_i$ with 
$\varphi^* H \cdot l_1=1$ and
$\varphi^* H \cdot l_i=0$ for $i \geq 2$. 

Since $\varphi^* H \cdot l_1=1$, the curve
$l_1$ is a subsection of a line, so it does not pass through $z$.
Thus there exists a $j \geq 2$ such that $z \in l_j$.
Since $\varphi^* H \cdot l_j=0$, this implies 
$l_j \subset \fibre{\varphi}{\varphi(z)}=\fibre{\varphi}{y}$.
We have shown that for a general point $y \in \Delta$, 
there exists a component of a degeneration of $l$
contained in $\fibre{\varphi}{y}$. 
Therefore there exists at least one irreducible
component of \fibre{\varphi}{\Delta} that is dominated 
by irreducible components of degenerations of $l$. 
Lemma \ref{lemmabirationalevaluation} shows that this is a contradiction. 
\finpreuve

\medskip

The following examples show the importance of the 
assumption that $e(X,l)=1$ for the theorem.

\medskip

\begin{examples} 
\label{examplespardini}
\begin{enumerate}
\item[(1)]
Let $X$ be a double cover of $\PP^1 \times \PP^2$
whose branch locus is a divisor of bidegree $(2,2)$. The threefold $X$ is Fano
with Picard number $2$, 
denote by $\psi~: X \to \PP^1$ and $\varphi~: X \to \PP^2$
the two natural projections. The map $\varphi$ is a conic bundle, whose
dicriminant is a quartic curve 
in $\PP^2$ and the map $\psi$ is a quadric bundle.
Let $l$ be a general line in $\PP^2$ and set $S_l := \varphi ^{-1}(l)$.
The surface $S_l$ is a del Pezzo surface  
and the induced map $\varphi~: S_l \to l \simeq \PP^1$
has exactly $4$ singular fibres. More precisely, we have a diagram

\begin{equation}
 \xymatrix{ 
S_l \ar[rd]^{\varphi} \ar[r]^{\mu} & {\mathbb F}_1 \ar[d] 
\ar[r]^{\mu_ 0} & \PP^2 \\
 & l\simeq \PP^1 & } 
\end{equation}

\noindent where $\mu_0$ is the blow-up of a point $x_0\in \PP^2$
and $\mu$ is the blow-up of $4$ points $p_1,\ldots,p_4$ in different fibres.
Set $\pi~:= \mu_0 \circ \mu~: S_l \to \PP^2$, and let $d$ be  
a general line in $\PP^2$ and $\omega = \pi^{-1}(d)$. 
Clearly $\omega$ is a quasi-line of $S_l$ and a section
of $\varphi~: S_l \to l$, therefore $\omega$ is a quasi-line of $X$ by lemma \ref{lemmafibrations}.  
One sees that $(X,\omega)$ satisfies all hypothesis of theorem \ref{theoremnotconicbundle} 
{\em except} $e(X,\omega) = 1$: if $x$ and $x'$ are two general points
of $\omega$, take a cubic curve in $P$, passing through 
$x$, $x'$, $\mu_0(p_i)$ for $i=1,\ldots,4$ with multiplicity $1$,
and through $x_0$ with multiplicity $2$. 
Then the strict transform of this cubic is also a quasi-line 
of $S_l$ and a section
of $\varphi~: S_l \to l$, hence is also a quasi-line of $X$ by 
lemma \ref{lemmafibrations}.  
  
Note that these two quasi-lines do not belong to 
the same family as quasi-lines of $S_l$ but belong to the same family as
quasi-lines of $X$!

\item[(2)] Exactly in the same manner, 
let $X$ be a double cover of $\PP^1 \times \PP^2$
whose branch locus is a divisor of bidegree $(2,4)$. The threefold $X$ is Fano
with Picard number $2$, 
denote by $\psi~: X \to \PP^1$ and $\varphi~: X \to \PP^2$
the two natural projections. The map $\varphi$ is a conic bundle, whose
discriminant is a curve of degree $8$ 
in $\PP^2$ and the map $\psi$ is a del Pezzo fibration.
Let $l$ be a general line in $\PP^2$ and $S_l := \varphi ^{-1}(l)$.
The surface $S_l$ is a surface with nef anti-canonical bundle
(isomorphic to a $\PP^2$ blown-up at $9$ general points) 
and the induced map $\varphi~: S_l \to l \simeq \PP^1$
has exactly $8$ singular fibres. One can show that $S_l$ contains 
a quasi-line, section
of $\varphi~: S_l \to l$, hence a quasi-line of $X$.   
\end{enumerate}
\end{examples}

\section{Fano bundles}
\label{subsectionfanobundles}

Theorem \ref{theoremnotconicbundle} shows that projective bundles over $\PP^m$ play a special
role in the classification of models $(X,l)$ with $e(X,l)=1$, in particular if $X$ is minimal with
respect to $l$. More generally we might ask when a Fano manifold that is a projectivised 
bundle $\PP(E)$ (a so-called
Fano bundle) over $\PP^m$ contains a quasi-line. 
Classification results and vector bundle techniques 
will allow us to give a complete answer if $\rk E=2$.
Together with theorem \ref{theoremnotconicbundle}, this also ends the proof of theorem~C.

\medskip

{\bf Notation.}
Let \holom{\varphi}{\PP(E)}{\PP^m} be a projectivised bundle over $\PP^m$
(we use here Grothendieck's definition: $\PP(E)$
is the variety of hyperplanes of $E$). {\em Throughout the whole 
section}, we will denote by $H$ a hyperplane divisor in $\PP^m$, and by 
$\xi_E$ the tautological divisor on $\PP(E)$. If $l \subset \PP(E)$ is a quasi-line, we will use
frequently $\varphi^* H \cdot l \geq 1$.

\begin{prop}
\label{propositionstability}
\begin{enumerate}
\item Let $(X,l)$ be a smooth model, and suppose that $X$ is minimal with respect to $l$.  
Suppose furthermore that $X \simeq \PP(E)$ for some rank 2 vector bundle   
$E$ over the projective space $\PP^m$. Then the vector bundle $E$ is stable. 
\item Let $E$ be a semistable vector bundle of rank 2 
with odd first Chern class over the projective space $\PP^m$. 
Then $X \simeq \PP(E)$ contains a quasi-line $l$ and
the natural projection
\holom{\varphi}{X}{\PP^m} is a morphism of models $(X,l) \rightarrow (\PP^m,line)$.
\end{enumerate}
\end{prop}

{\em Proof.}
For the first statement, we
argue by contradiction and suppose that $E$ is not stable. 
Then, up to twisting $E$ with a line bundle, we can suppose without loss of generality  
that $-1 \leq c_1(E) \leq 0$ and $h^0(X,E)>0$ (cf. \cite[Rem.3.0.1]{Ha78}).
Let \holom{\varphi}{X}{\PP^m} be the projection map,  then
\[
- K_X = (m+1-c_1(E)) \varphi^* H + 2 \xi_E.
\]
It follows that
\begin{equation}
\label{equationcomputation}
m + 2 = - K_X \cdot l = 
(m+1-c_1(E)) \varphi^* H \cdot l + 2 \xi_E \cdot l
\geq  m+1 + 2 \xi_E \cdot l. 
\end{equation}
Since $h^0(X, \xi_E)=h^0(\PP^m,E)>0$ and $X$ is minimal, we have $\xi_E \cdot l \geq 1$. 
This implies
\[
m + 2 \geq m+1 + 2 \xi_E \cdot l \geq m + 1 + 2,
\]
a contradiction. 

For the second statement, let $l' \subset \PP^m$ be a general line. 
By the Grauert-M\"ulich theorem \cite[p.206]{OSS80}, we have
$E|_{l'} \simeq \sO_{\PP^1}(a_1) \oplus \sO_{\PP^1}(a_2)$ with $|a_2 - a_1| \leq 1$.
Since $E$ has an odd first Chern class, we can suppose up to renumbering that $a_2=a_1+1$.
By \cite[Prop. 4.2]{IV03} there exists a quasi-line $l$ on $X$ such that $\holom{\varphi|_l}{l}{l'}$
is an isomorphism.
\finpreuve

\medskip

{\bf Remark.} By a well-known conjecture of Hartshorne \cite{Ha74}, there are no stable rank 2 vector bundles
over $\PP^m$ if $m \geq 7$. This shows (at least conjecturally) that the minimality with respect to a
quasi-line is rather restrictive.

\begin{prop}
\label{propositiondirectsums}
Let $E \simeq \oplus_{i=1}^{r+1} \sO_{\PP^m}(a_i)$ be a sum of line bundles on $\PP^m$ such that
$0 \leq a_1 \leq a_2 \leq \ldots \leq a_{r+1}$. Set $X := \PP(E)$, then $X$  
contains a quasi-line $l$ if and only if up to twisting
\[
E \simeq \sO_{\PP^m}^{\oplus r} \oplus \sO_{\PP^m}(1).
\]
In this case the natural map
\holom{\varphi}{X}{\PP^m} is a morphism of models $(X,l) \rightarrow (\PP^m,line)$.
\end{prop}

{\em Proof.}
The if part is immediate from \cite[Prop.4.2]{IV03}, for the only if part
we follow closely the proof from \cite[Prop.4.1]{IV03}. 

We have
\[
-K_X = (m+1-\sum_{i=1}^{r+1} a_i) \varphi^* H + (r+1) \xi_E.
\]
and
\[
h^0(X, \sO_X(\xi_E-a_{r+1} \varphi^* H)) = h^0 (\PP^m, E \otimes \sO_{\PP^m}(-a_{r+1}))>0,
\]
so the divisor $\xi_E-a_{r+1} \varphi^* H$ is linearly equivalent to an effective divisor. In particular
\[
(\xi_E-a_{r+1} \varphi^* H) \cdot l \geq 0.
\]
We want to show that even $(\xi_E-a_{r+1} \varphi^* H) \cdot l = 0$ and 
argue by contradiction.  Then 
\[
(\xi_E-a_{i} \varphi^* H) \cdot l \geq (\xi_E-a_{r+1} \varphi^* H) \cdot l \geq 1
\]
for all $i=1, \ldots r+1$, so
\begin{eqnarray*}
m+r+1 = -K_X \cdot l &=& (m+1-\sum_{i=1}^{r+1} a_i) \varphi^* H \cdot l + (r+1) \xi_E \cdot l 
\\
&=&
(m+1) \varphi^* H \cdot l + \sum_{i=1}^{r+1} (\xi_E-a_i \varphi^* H) \cdot l 
\\
& \geq & (m+1) + (r+1) 
\end{eqnarray*}
yields a contradiction.
So $(\xi_E-a_{r+1} \varphi^* H) \cdot l = 0$ and 
\[
h^0(X, \sO(\xi_E-a_{r+1} \varphi^* H)) = h^0 (\PP^m, E \otimes \sO(-a_{r+1}))=1,
\]
since a quasi-line has strictly positive intersection
number with every divisor that moves. This implies that $a_i < a_{r+1}$ for $i=1, \ldots, r$, so
\[
(\xi_E-a_{i} \varphi^* H) \cdot l > (\xi_E-a_{r+1} \varphi^* H) \cdot l = 0
\]
for all $i=1, \ldots r$.
We repeat the preceeding computation
\begin{eqnarray*}
m+r+1 = -K_X \cdot l
&=&
(m+1) \varphi^* H \cdot l + \sum_{i=1}^{r+1} (\xi_E-a_i \varphi^* H) \cdot l 
\\
& \geq & (m+1) \varphi^* H \cdot l + r, 
\end{eqnarray*}
so $\varphi^* H \cdot l=1$. This implies that $\varphi(l)$ is a line in $\PP^m$ and  
\holom{\varphi|_l}{l}{\varphi(l)} is an isomorphism. We conclude with a second application of \cite[Prop.4.2]{IV03}. \finpreuve

\begin{cor}
\label{corollaryfanobundlem3}
Let $(X,l)$ be a smooth Fano model. 
Suppose that $X \simeq \PP(E)$ for some rank 2 vector bundle   
$E$ over the projective space $\PP^m$ with $m \geq 3$. 
Then up to twisting with a line bundle we have 
\[
E \simeq \sO_{\PP^m}(-1) \oplus \sO_{\PP^m},  
\]
so $X$ is the blow-up of $\PP^{m+1}$ in a point and $l$ is the preimage of a line in $\PP^{m+1}$. In particular $X$
is not minimal. 
\end{cor}

{\em Proof.} 

{\it 1st case: $m\geq 4$.}
Since $E$ is a rank 2 Fano bundle, it splits by 
\cite[Main thm.]{APW94} in a direct sum of line bundles. 
Conclude with proposition \ref{propositiondirectsums}. 

{\it 2nd case: $m=3$.}
We normalise $E$ such that $c_1(E)=0$ or $c_1(E)=-1$.
If the first Chern class is zero, we have
\[
- K_X = 4 \varphi^* H +  2 \xi_E.
\]
In particular $X$ has index 2, so it is clear that $X$ does not contain a quasi-line ($-K_X \cdot l=5$).
If $c_1(E)=-1$, 
the classification of Fano bundles of rank 2 on $\PP^3$ \cite[Thm.2.1]{SW90b},
implies that $E \simeq \sO_{\PP^3} \oplus \sO_{\PP^3}(-1)$ or
$E \simeq \sO_{\PP^3} (-2) \oplus \sO_{\PP^3}(1)$. 
The second case is excluded by proposition \ref{propositiondirectsums}. \finpreuve

\begin{prop}
\label{propositionp2fanobundle}
Let $(X,l)$ be a smooth Fano model. 
Suppose that $X \simeq \PP(E)$ for some vector bundle    
$E$ of rank $r$ over the projective space $\PP^2$. 
Then up to twisting with a line bundle we have either
\[
E \simeq \sO_{\PP^2}^{\oplus r-1} \oplus \sO_{\PP^2}(1),
\]
or $E$ is defined by an exact sequence
\[
0 \rightarrow \sO_{\PP^2}(-1)^{\oplus 2} \rightarrow \sO_{\PP^2}^{\oplus r+2} \rightarrow E \rightarrow 0,
\]
or
\[
E \simeq T_{\PP^2}(-1) \oplus \sO_{\PP^2}^{\oplus r-2},
\]
or $E$ is defined by an exact sequence
\[
0 \rightarrow \sO_{\PP^2}(-2) \rightarrow \sO_{\PP^2}^{r+1} \rightarrow E \rightarrow 0.
\]
Vice versa, if $E$ is a vector bundle over $\PP^2$ as above, then $\PP(E)$ is Fano and contains a quasi-line. 

In the first two cases 
there exists a birational morphism $X \rightarrow \PP^{r+1}$ that induces an equivalence of
models $(X,l) \simeq (\PP^{r+1}, line)$; in the last two cases $X$ is minimal with respect to $l$.
In the third case $e(X,l)=1$ and in the fourth case $e(X,l)>1$.
\end{prop}

{\em Proof.} 
Denote by \holom{\varphi}{X}{\PP^2} the projection map.
We normalise the vector bundle $E$ such that $0 \leq c_1(E) \leq r - 1$, then
by \cite[1.6.]{SW90} the Fano condition implies $0 \leq c_1(E) \leq 2$. Furthermore
the vector bundle $E$ is globally generated \cite[Prop.2.1]{SW90}. Since $E$ can't be trivial
(the product $\PP^{r-1} \times \PP^2 \simeq \PP(E)$ does not contain a quasi-line), this implies
\[
(*) \qquad h^0(X, \sO_{X}(\xi_E)) = h^0(\PP^2, E) \geq r+1,
\]
in particular
\[
\xi_E \cdot l \geq 1.
\]
Therefore
$\varphi^* \sO_{\PP^2}(3-c_1(E)) \cdot l \geq 1$ implies
\begin{eqnarray*}
r+2 &=& -K_X \cdot l 
\\
&=& 
(3-c_1(E)) \varphi^* H \cdot l + r \xi_E \cdot l
\\
&\geq & 1 + r \xi_E \cdot l.
\end{eqnarray*}
Hence $\xi_E \cdot l=1$, since $\sO_X(\xi_E)$ is globally generated it induces
by \cite[Thm.1.12]{IV03} a surjective morphism with connected fibres 
\holom{\psi}{X}{\PP^N:=\PP(H^0(X,\sO_X(\xi_E)))} 
that induces a morphism of models $(X,l) \rightarrow (\PP^N,line)$.
In particular
\[
(**) \qquad r + 2 = \dim X+1 \geq h^0(X, \sO_{\PP(E)}(1)) = h^0(\PP^2, E)
\]
and equality holds if and only if 
$\psi$ induces an equivalence of models $(X,l) \simeq (\PP^{r+1},line)$.
In view of the inequalities $(*)$ and $(**)$ we have to treat two cases.
 
{\it The case $h^0(\PP^2, E)=r+2$.} 
Going through the list in \cite[Thm.]{SW90} yields that in this case, we have
$E \simeq \sO_{\PP^2}^{\oplus r-1} \oplus \sO_{\PP^2}(1)$ or $E$ is defined by an exact sequence
\[
0 \rightarrow \sO_{\PP^2}(-1)^{\oplus 2} \rightarrow \sO_{\PP^2}^{\oplus r+2} \rightarrow E \rightarrow 0.
\]
In both cases, the manifold $X$ is a blow-up of $\PP^{r+1}$, so it is clear that $X$ contains
a quasi-line and is not minimal with respect to it. 

{\it The case $h^0(\PP^2, E)=r+1$.}
In this case $\psi$ gives 
a fibration of relative dimension one onto $\PP^r$ that 
induces a morphism of models $(X,l) \rightarrow (\PP^r,line)$.
Using the list in \cite[Thm.]{SW90}, we then see that 
$E \simeq T_{\PP^2}(-1) \oplus \sO_{\PP^2}^{\oplus r-1}$
or $E$ is defined by an exact sequence
\[
0 \rightarrow \sO_{\PP^2}(-2) \rightarrow \sO_{\PP^2}^{\oplus r+1} \rightarrow E \rightarrow 0.
\]
In the first case $c_1(E)=1$, so 
\[
r+2 = - K_X \cdot l = 2 \varphi^* H \cdot l + r.
\]
implies that $\varphi^* H \cdot l = 1$. Hence $\varphi(l)$ is a line in $\PP^2$ and $l$ is a section.
Vice versa \cite[Prop.4.3]{IV03} shows that $\PP(T_{\PP^2}(-1) \oplus \sO_{\PP^2}^{\oplus r-1})$
contains a quasi-line $l$ that is a subsection of a line in $\PP^2$,
Since the two elementary contractions are of fibre type, the manifold $X$ 
is minimal with respect to $l$ by lemma \ref{lemmabirationalcontractions}, and
$e(X,l)=1$ follows from the remark after the proof.

In the second case $X$ can be realised as a smooth divisor of degree (2,1) in $\PP^2 \times \PP^r$ \cite[Thm.]{SW90},
and $\varphi$ (resp. $\psi$) identifies to the projection on the first (resp. second) factor. 
Since the property of containing a quasi-line is stable under small deformations \cite[Prop.3.10]{BBI00}, 
we can choose $X$ sufficiently
general such that $\psi$ has no higher-dimensional fibres, {\em i.e.} 
yields a conic bundle structure on $X$
and a straightforward computation shows that the discriminant locus has degree 3.
By adjunction $\omega_X^* \simeq  \sO_X(1,r)$, so if $l' \subset \PP^{r}$ is a general
line and $S:=\fibre{\psi}{l'}$ is
its preimage, repeated adjunction yields $\omega_S^* \simeq \sO_S(1,1)$. Hence $S$ is a del Pezzo surface
and \holom{\psi|_S}{S}{l'} has three singular fibres. 
The surface $S$ is a blow-up of $\PP^2$ in four points, so a general line in $\PP^2$ yields a quasi-line $l$
in $S$ that is a section of $\psi|_S$. 
By lemma \ref{lemmafibrations} this shows that $l$ is a quasi-line in $X$. 
Since the two elementary contractions are of fibre type, the manifold $X$ 
is minimal with respect to $l$ by lemma \ref{lemmabirationalcontractions}.
Since $\psi$ is not a smooth morphism, theorem \ref{theoremnotconicbundle} shows that $e(X,l)>1$. \finpreuve

\medskip

{\bf Remark.}
Let \holom{\varphi}{(X,l)}{(Y,\varphi(l))} be a fibration between smooth models, and choose a sufficiently
general deformation of $l$ such that $Z:=\fibre{\varphi}{\varphi(l)}$ is smooth. We define
\[
e(Z) := \sum_{l \in \mathcal H} e(X,l),
\]
where the sum goes over all the irreducible components of the Chow scheme such that the general
point parametrises a quasi-line. Then it is not hard to see 
(see \cite[Lemma 5.9]{IN03} for the special case of
a projective bundle) that
\[
e(X,l) \leq e(Z) \cdot e(Y, \varphi(l)).
\]
In particular if $X \simeq \PP(E)$ for a vector bundle $E$ of rank $r$
over $Y$, then $Z \simeq \PP(\sO_{\PP_1}^{\oplus r-1} \oplus \sO_{\PP_1}(1))$ \cite[Prop. 4.1, Prop.4.2]{IV03}, 
so $e(Z)=1$. 

The fourth case of the proposition shows that not all the Fano bundles over $\PP^2$ satisfy $e(X,l)=1$.
In particular $X \rightarrow \PP^2$ is not a morphism of models! 
Note furthermore that 
the manifold $\PP(T_{\PP^2} (-1) \oplus \sO_{\PP}^{\oplus r-2})$ contains 
a quasi-line and is minimal, but for
$r \geq 3$, the vector bundle  
$T_{\PP^2} (-1) \oplus \sO_{\PP}^{\oplus r-2}$ isn't even semistable.
This shows that for a generalisation of proposition \ref{propositionstability} to bundles
of higher rank new ideas are necesssary.   

\medskip

{\bf Towards a classification of Fano models ?} It is a natural question
to ask for a complete list of Fano models $(X,l)$ such that $X$ is minimal with respect to $l$.
This turns out to be a rather lengthy exercice, in fact using the techniques in this paper
(in particular the examples \ref{examplespardini}) one can show that out of the nine
primitive Fano threefolds with Picard number two, all but $\PP^2 \times \PP^1$ and
$\PP(\sO_{\PP^2} \oplus \sO_{\PP^2}(2))$ contain a quasi-line. Among these models, five 
have only fibre type contractions, so they are minimal with respect to the quasi-line 
(lemma \ref{lemmabirationalcontractions}).
Subsection \ref{subsectionexamples} shows how to obtain more examples in the non-primitive case.
This convinces us that classification efforts should concentrate on the case $e(X,l)=1$.
In this case we only know one example with Picard number at least two: the flag manifold $\PP(T_{\PP^2})$.

{\em \small
\noindent Laurent Bonavero. 
Institut Fourier, UMR 5582, 
Universit\'e de Grenoble 1, BP 74. 
38402 Saint Martin d'H\`eres. France. 
\\
e-mail : laurent.bonavero@ujf-grenoble.fr

\noindent Andreas H\"oring. 
IRMA, Universit\'e Louis Pasteur Strasbourg, 7 rue Ren\'e Descartes. 
67084 Strasbourg. France\\
\noindent e-mail : ahoering@ujf-grenoble.fr
}

\end{document}